\documentclass[12pt]{article}
\usepackage{amsmath,amssymb}
\usepackage[mathscr]{euscript}
\usepackage[OT1,T1]{fontenc}
\usepackage[all]{xy}
\usepackage{hyperref}
\usepackage{cancel}
\usepackage[normalem]{ulem}

\usepackage{xcolor}

\newtheorem{theo}{Theorem}[section]

\newcommand{\dR}{\hbox{d}}
\newcommand{\goG}{\mathfrak{G}}
\newcommand{\gog}{\mathfrak{g}}
\newcommand{\gos}{\mathfrak{s}}
\newcommand{\gol}{\mathfrak{l}}
\newcommand{\gou}{\mathfrak{u}}

\newcommand{\celg}{{\underline{\mathfrak{g}}}}
\newcommand{\cell}{{\underline{\,\mathfrak{l}\,}}}
\newcommand{\cels}{{\underline{\mathfrak{s}}}}
\newcommand{\celu}{{\underline{\mathfrak{u}}}}

\title{A note on a Kaluza--Klein theory}
\author{Fr{\'e}d{\'e}ric \textsc{H{\'e}lein}\footnote{Institut de Math{\'e}matiques de Jussieu--Paris Rive Gauche,
UMR CNRS 7586, Universit{\'e} Paris Cité,
UFR de Math{\'e}matiques,  B{\^a}timent Sophie Germain,
8 place Aurélie Nemours, 75013 Paris, France, \textsf{helein@math.univ-paris-diderot.fr}}
}

\begin{document}

\maketitle
\textbf{Abstract} --- \emph{We present a variant of the Kaluza-Klein theory that was introduced in previous papers  (\cite{helein2020,helein22,helein25}), based on an action functional whose critical points, in the presence of point-wise constraints, lead—under reasonable assumptions—to the existence of a principal fiber bundle and solutions to the Einstein-Yang-Mills system. The difference is that it is not necessary to impose a point-wise constraint.}\\

In \cite{helein2020,helein22,helein25} we introduced
a Kaluza--Klein type theory for gravity coupled with gauge fields which does not require to assume \emph{a priori} the existence of a principal bundle structure.
In order to achieve such a theory on an $n$-dimensional spacetime with gauge fields with coefficients in an $r$-dimensional Lie algebra $\gog$, the basic ingredients are an $N$-dimensional manifold $\mathcal{Y}$, where $N = n+r$, and fields defined on it which consist in a collection of exterior differential forms of degrees 1 and $N-2$, \emph{satisfying some pointwise constraints} (see (\ref{constrainttoberemoved}) below). An action functional $\mathscr{A}$ is defined on the set of such fields and it is shown in, e.g., \cite{helein25} that, under reasonable hypotheses, each critical point of this functional provides us with:
\begin{enumerate}
 \item a foliation of $\mathcal{Y}$ by $r$-dimensional leaves
\end{enumerate}
Morever, if $\gog$ is the Lie algebra of a compact simply connected Lie group $\widehat{\goG}$:
\begin{enumerate}
 \item[2.]  the set of leaves of the foliation has the structure of an $n$-dimension manifold $\mathcal{X}$, which can be identified with a model for the space-time.
 \item[3.]  $\mathcal{Y}$ acquires the structure of a principal bundle over $\mathcal{X}$ (the fibers being the leaves of ths foliation)
 \item[4.] a metric and a $\gog$-valued gauge field on $\mathcal{X}$ which are solutions of the Einstein--Yang--Mills system of equations in presence of an explicitely given cosmological constant.
\end{enumerate}
We remark that the hypothesis that $\gog$ be the Lie algebra of a compact simply connected Lie group is satisfied by $SU(2)$, $SU(3)$ and $SU(2)\times SU(3)$ but not by $U(1)$. However similar conclusions are also valid for $U(1)$ under additional hypotheses (see \cite{helein22,helein25}).

The result in this note is that the pointwise constraint (\ref{constrainttoberemoved}) can be removed by adding an extra field playing the role of a Lagrange multiplier. We will also observe that the resulting action functional is the sum of a Palatini type and of a BF term, in which the Lie bracket is replaced by a bracket whose structure coefficients are themselves partially dynamical. Also the new model presented in this note simplifies the computation of the Euler--Lagrange equations.

\section{Our previous Kaluza--Klein model}
In the following we use the same notations and conventions as in \cite{helein25}.

Let $\gos$ be an $n$-dimensional vector space which is a model for a flat space-time (e.g. the Minkowski space).
Let us endow $\gou = \gos\oplus \gog$ with the Lie bracket $[\cdot,\cdot]$ which extends trivially the Lie algebra structure of $\gog$ (i.e such that $\gos$ belongs to the center of $(\gou,[\cdot,\cdot])$). Assume that $\gou$ is endowed with a symmetric bilinear form $\textsf{h} = \textsf{b}\oplus \textsf{k}$ which is invariant by the adjoint action of $\gog$ on $\gou$ and such that $\gos\perp \gog$. We denote by $\gol$ the Lie algebra of the group $SO(\gou,\textsf{h})$ of special linear isometries of $(\gou, \textsf{h})$.

Let $\mathcal{Y}$ be a connected, oriented manifold of dimension $N:= n+r$.
Consider the space of fields\[
\begin{array}{r}
\mathscr{E}:=
 \{\, (\theta^\gou,\varphi^\gol,\pi_\gou)\in \left(\gou \otimes\Omega^1(\mathcal{Y})\right)\times \left(\gol \otimes\Omega^1(\mathcal{Y})\right) \times \left(\gou^*\otimes\Omega^{N-2}(\mathcal{Y})\right)\,; \\
 \hbox{rank}\,\theta^\gou = N,
(\theta^\gou,\varphi^\gol,\pi_\gou)\hbox{ are of class }\mathscr{C}^2\hbox{ and satisfy }(\ref{constrainttoberemoved}) \}
\end{array}
\]
where
\begin{equation}\label{constrainttoberemoved}
 \theta^\mathfrak{s}\wedge \theta^\mathfrak{s} \wedge \pi_\gou =0
\end{equation}
Let $\kappa{_\gol}^{\gou\gou}\in \gol^*\otimes\gou\otimes\gou$
be defined as by (16) in \cite{helein25}
% (\ref{kappasimple})
(this tensor is invariant by $\hbox{Ad}_{\widehat{\goG}}$) and set
$\Phi^\gou:= \dR\theta^\gou + \frac{1}{2}[\theta^\gou\wedge \theta^\gou]$,
$\Phi^\gol:= \dR\varphi^\gol+
 \frac{1}{2}[\varphi^\gol\wedge \varphi^\gol]$ and
 $\Phi^{\gou\gou}:= \kappa{_\cell}^{\gou\gou}\Phi^\cell$.
Here and in the following the underlined indices $\cels,\celu,\celg$ indicate summations over repeated indices, see \cite{helein25}, \S 3.3 for details. Let $\theta^{(N)}:= \theta^1\wedge \cdots \wedge \theta^N$, where $(\theta^1, \cdots, \theta^N)$ is a decomposition of $\theta^\gou$ in a basis of $\gou$. Let $\Lambda_0$ be some real constant.
Consider the action functional $\mathscr{A}$ defined on  $\mathscr{E}$ by:
\begin{equation}\label{actionKK}
 \mathscr{A}[\theta^\gou,\varphi^\gol,\pi_\gou]:= \int_{\mathcal{Y}}  \pi_\celu \wedge \Theta^\celu + \frac{1}{2}\theta_{\celu\celu}^{(N-2)} \wedge \Phi^{\celu\celu} - \Lambda_0 \theta^{(N)}
\end{equation}
Then, if $\widehat{\goG}$ is compact, any critical point $(\theta^\gou,\varphi^\gol,\pi_\gou) \in \mathscr{E}
$ satisfies the following system of equations
\begin{equation}\label{ELFrobeniusdetaillee}
 \left\{
 \begin{array}{ccl}
\hbox{d}\theta^\gos & = &  \frac{1}{2}\Theta{^\gos}_{\cels\cels}\,\theta^\cels\wedge \theta^\cels \\
\hbox{d}\theta^\gog +\frac{1}{2}[\theta^\gog\wedge \theta^\gog] & = &  \frac{1}{2}\Theta{^\gog}_{\cels\cels}\,\theta^\cels\wedge \theta^\cels
 \end{array}\right.
\end{equation}
\begin{equation}\label{ELsimple}
 \hbox{d}^\varphi \theta^\gou := \hbox{d}\theta^\gou + (\rho(\varphi)\wedge \theta)^\gou = 0
\end{equation}
where $\rho$ denotes the standard representation of $\gol$ on $\gou$, and
\begin{equation}\label{ELEinsteinbrut}
 \dR^\theta\pi_\gou +
 \frac{1}{2}\theta^{(N-3)}_{\celu\celu\gou}\wedge \Phi^{\celu\celu}
- \Lambda_0 \theta^{(N-1)}_\gou
= \Theta{^\celu}_{\gou\cels}\pi{_\celu}^{\celg\cels}\ \theta^{(N-1)}_{\celg}
\end{equation}
As a consequence of these equations we can associate an $n$-dimensional manifold $\mathcal{X}$ and a submersion $\pi:\mathcal{Y}\longrightarrow \mathcal{X}
$ such that $\mathcal{Y}$ has the structure of a $\goG$-principal bundle over $\mathcal{X}$, endowed with a connection and a pseudo-Riemannian metric on $\mathcal{X}$ which are solutions of the Einstein--Yang--Mills system of equations \cite{helein2020, helein22,helein25}.

\section{Removing the constraint}

The constraint (\ref{constrainttoberemoved}) in the definition of $\mathscr{E}$ can be removed if we add an extra field $\lambda{^\gou}_{\gos\gos}\in  \gou\otimes \gos^*\wedge \gos^*\otimes \Omega^0(\mathcal{Y})$, playing the role of a Lagrange multiplier.
The space of fields is then:
\[
\begin{array}{ll}
\widetilde{\mathscr{E}}:=
 \{\, (\theta^\gou,\varphi^\gol,\pi_\gou,\lambda{^\gou}_{\gos\gos}) \in \left(\gou \otimes\Omega^1(\mathcal{Y})\right)\times \left(\gol \otimes\Omega^1(\mathcal{Y})\right)\times \left(\gou^*\otimes\Omega^{N-2}(\mathcal{Y})\right) \times \left(\gou\otimes \gos^*\wedge \gos^*\otimes \Omega^0(\mathcal{Y})\right)\\
\quad\quad\quad \quad \hbox{rank}\,\theta^\gou = N,\hbox{ all fields }\mathscr{C}^2 \}
\end{array}
\]
We then define on $\widetilde{\mathscr{E}}$ the action functional $\widetilde{\mathscr{A}}$ by:
\begin{equation}\label{actionKKplus}
 \widetilde{\mathscr{A}} =[\theta^\gou,\varphi^\gol,\pi_\gou, \lambda{^\gou}_{\gos\gos}]:= \int_{\mathcal{Y}}  \pi_\celu \wedge \left( \Theta^\celu + \frac{1}{2}\lambda{^\celu}_{\cels\cels}\theta^\cels\wedge \theta^\cels\right) + \frac{1}{2}\theta_{\celu\celu}^{(N-2)} \wedge \Phi^{\celu\celu} - \Lambda_0 \theta^{(N)}
\end{equation}
Alternatively we may define the $\lambda$-depending skewsymmetric bracket $[\cdot,\cdot]_\lambda:\gou\times \gou \longrightarrow \gou$ by
\begin{equation}
 \forall \xi^\gou,\eta^\gou\in \gou,\quad [\xi^\gou,\eta^\gou]_\lambda:= \lambda{^\gou}_{\cels\cels}\xi^\cels \eta^\cels + \mathbf{c}{^\gog}_{\celg\celg} \xi^\celg \eta^\celg
\end{equation}
and set
\[
 [\theta^\gou\wedge \theta^\gou]_\lambda = \lambda{^\gou}_{\cels\cels}\theta^\cels\wedge \theta^\cels + \mathbf{c}{^\gog}_{\celg\celg} \theta^\celg\wedge \theta^\celg
\]
so that
\[
  \widetilde{\mathscr{A}} [\theta^\gou,\varphi^\gol,\pi_\gou, \lambda{^\gou}_{\gos\gos}]:= \int_{\mathcal{Y}}  \pi_\celu \wedge \left( \dR\theta^\celu + \frac{1}{2}[\theta^\gou\wedge \theta^\gou]^\celu_\lambda\right) + \frac{1}{2}\theta_{\celu\celu}^{(N-2)} \wedge \Phi^{\celu\celu} - \Lambda_0 \theta^{(N)}
\]
\begin{theo}
Let $(\theta^\gou,\varphi^\gol,\pi_\gou,\lambda{^\gou}_{\gos\gos}) \in \widetilde{\mathscr{E}}$.
Then $(\theta^\gou,\varphi^\gol,\pi_\gou,\lambda{^\gou}_{\gos\gos})$ is a critical point of $\widetilde{\mathscr{A}}$ if and only if $(\theta^\gou,\varphi^\gol,\pi_\gou) \in \mathscr{E}$ is a critical point of $\mathscr{A}$ under the constraint (\ref{constrainttoberemoved}) and if
\[
 \Theta{^\gou}_{\gos\gos} + \lambda{^\gou}_{\gos\gos} = 0
\]
\end{theo}
\emph{Proof} ---
For simplicity we assume that $\Lambda_0 = 0$. We derive the Euler--Lagrange equations of $\widetilde{\mathscr{A}}$ and show that we recover the  Euler--Lagrange equations of $\mathscr{A}$ given in \cite{helein25}.\\
\noindent
\textbf{First variation with respect to coefficients of $\pi_\gou$} --- The equation is much simpler to derive: the action functional is stationary with respect to first order variations $(\theta^\gou,\varphi^\gol,\pi_\gou,\lambda{^\gou}_{\gos\gos})\longmapsto (\theta^\gou,\varphi^\gol,\pi_\gou + \varepsilon \delta\pi_\gou,\lambda{^\gou}_{\gos\gos})$ if and only if
\[
\dR\theta^\gou + \frac{1}{2}[\theta^\gou\wedge \theta^\gou]_\lambda = \dR\theta^\gou + \frac{1}{2}\lambda{^\gou}_{\cels\cels}\theta^\cels\wedge \theta^\cels + \frac{1}{2}\mathbf{c}{^\gog}_{\celg\celg} \theta^\celg\wedge \theta^\celg= 0
\]
This is equivalent to
\[
\left\{
 \begin{array}{ccl}
\hbox{d}\theta^\gos + \frac{1}{2}\lambda{^\gos}_{\cels\cels}\theta^\cels\wedge \theta^\cels & = &  0 \\
\hbox{d}\theta^\gog +\frac{1}{2}[\theta^\gog\wedge \theta^\gog] + \frac{1}{2}\lambda{^\gog}_{\cels\cels}\theta^\cels\wedge \theta^\cels   & = &  0
 \end{array}\right.
\]
or, by defining the coefficients $\Theta{^\gou}_{\gou\gou}$ by $\dR\theta^\gou + \frac{1}{2}[\theta^\gou\wedge \theta^\gou] = \frac{1}{2}\Theta{^\gou}_{\celu\celu}\theta^\celu\wedge \theta^\celu$,
\begin{equation}\label{equA}
\Theta{^\gou}_{\gos\gos} + \lambda{^\gou}_{\gos\gos} = 0
\quad \hbox{ and }\quad
\Theta{^\gou}_{\gos\gog} = \Theta{^\gou}_{\gog\gog} = 0
\end{equation}
which also implies (\ref{ELFrobeniusdetaillee}) by eliminating $\lambda{^\gou}_{\gos\gos}$.

\noindent
\textbf{First variation with respect to $\varphi^\gol$} --- The computation is the same as in \cite{helein25} and leads to the condition that $\dR^\varphi \theta^{(N-2)}_{\gou\gou} = 0$, which, for $N>2$, is equivalent to the equation
\begin{equation}\label{equB}
\dR^\varphi \theta^\gou = 0
\end{equation}
We thus recover (\ref{ELsimple}), which means that the connection on $T\mathcal{Y}$ defined by $\varphi^\gol$ is torsion free.\\

\noindent
\textbf{First variation with respect to $\theta^\gou$} --- Again the computation is simpler. We look at variations $(\theta^\gou,\varphi^\gol,\pi_\gou,\lambda{^\gou}_{\gos\gos})\longmapsto (\theta^\gou + \varepsilon\delta\theta^\gou, \varphi^\gol, \pi_\gou ,\lambda{^\gou}_{\gos\gos})$, for any $\delta\theta^\gou = \tau^\gou$ with compact support.
Through these variations,
\[
 \left\{
 \begin{array}{cclll}
  \theta_{\gou\gou}^{(N-2)} & \longmapsto &\theta_{\gou\gou}^{(N-2)} + \varepsilon\delta\theta_{\gou\gou}^{(N-2)} + o(\varepsilon) & \hbox{with} & \delta\theta_{\gou\gou}^{(N-2)} = \tau^\celu\wedge \theta^{(N-3)}_{\gou\gou\celu} \\
  \Theta^\gou & \longmapsto & \Theta^\gou + \varepsilon\delta\Theta^\gou+ o(\varepsilon) & \hbox{with} & \delta\Theta^\gou = \dR^\theta \tau^\gou = \dR \tau^\gou + [\theta^\gou,\tau^\tau]
 % \theta^{(N)} & \longmapsto & \theta^{(N)} + \varepsilon\delta\theta^{(N)}+ o(\varepsilon) & \hbox{with} & \delta\theta^{(N)} = \tau^\celu\wedge \theta^{(N-1)}_\celu
 \end{array}
\right.
\]
Thus
\[
\begin{array}{ccl}
  \delta\widetilde{\mathscr{A}} & = & \displaystyle \int_\mathcal{Y}  \pi_\celu \wedge \left( \dR^\theta\tau^\celu + \lambda{^\celu}_{\cels_1\cels_2}\tau^{\cels_1}\wedge \theta^{\cels_2}\right) + \frac{1}{2}\tau^\celu \wedge\theta_{\celu_1\celu_2\celu}^{(N-3)} \wedge \Phi^{\celu_1\celu_2} \\
  %%%%%
  & = & \displaystyle \int_\mathcal{Y} \dR(\tau^\celu\wedge \pi_\celu) + \tau^\celu\wedge \left( \dR^\theta\pi_\celu +  \frac{1}{2}\theta_{\celu_1\celu_2\celu}^{(N-3)} \wedge \Phi^{\celu_1\celu_2}\right) +  \tau^{\cels}\wedge \lambda{^\celu}_{\cels\cels_2}\theta^{\cels_2}\wedge \pi_\celu
\end{array}
\]
Thus $\delta\widetilde{\mathscr{A}}$ vanishes if and only if
\begin{equation}\label{equC}
 \left\{
 \begin{array}{ccl}
\dR^\theta\pi_\gos + \frac{1}{2}\theta_{\celu_1\celu_2\gos}^{(N-3)} \wedge \Phi^{\celu_1\celu_2} + \lambda{^\celu}_{\gos\cels}\theta^{\cels}\wedge \pi_\celu & = & 0 \\
\dR^\theta\pi_\gog + \frac{1}{2}\theta_{\celu_1\celu_2\gog}^{(N-3)} \wedge \Phi^{\celu_1\celu_2} & = & 0
 \end{array}\right.
\end{equation}
\noindent
\textbf{First variation with respect to $\lambda{^\gou}_{\gos\gos}$} --- It gives us straightforwardly the constraint (\ref{constrainttoberemoved}),
which is equivalent to
\[
 \pi_\gou = \pi{_\gou}^{\cels_1\celg_2}\theta^{(N-2)}_{\cels_1\celg_2} + \frac{1}{2}\pi{_\gou}^{\celg_1\celg_2}\theta^{(N-2)}_{\celg_1\celg_2}
\]

\noindent
\textbf{Recovering the Euler--Lagrange equations of $\mathscr{E}$}\\
We already recovered (\ref{ELFrobeniusdetaillee}) and (\ref{ELsimple}).
By plugging $\Theta{^\gou}_{\gos\gos} + \lambda{^\gou}_{\gos\gos} = 0$ from (\ref{equA}) into (\ref{equC}) we get
\begin{equation}\label{equC2}
 \left\{
 \begin{array}{ccl}
\dR^\theta\pi_\gos + \frac{1}{2}\theta_{\celu_1\celu_2\gos}^{(N-3)} \wedge \Phi^{\celu_1\celu_2} & = & \Theta{^\celu}_{\gos\cels}\theta^{\cels}\wedge \pi_\celu \\
\dR^\theta\pi_\gog + \frac{1}{2}\theta_{\celu_1\celu_2\gog}^{(N-3)} \wedge \Phi^{\celu_1\celu_2} & = & 0
 \end{array}\right.
\end{equation}
However we can compute the right hand side of (\ref{equC2}) by using (\ref{constrainttoberemoved}):
\[
 \Theta{^\celu}_{\gos\cels}\theta^{\cels}\wedge \pi_\celu = \Theta{^\celu}_{\gos\cels}\theta^{\cels}\wedge (\pi{_\celu}^{\cels_1\celg_2}\,\theta^{(N-2)}_{\cels_1\celg_2} + \frac{1}{2}\pi{_\celu}^{\celg_1\celg_2}\,\theta^{(N-2)}_{\celg_1\celg_2})
 = - \Theta{^\celu}_{\gos\cels}\pi{_\celu}^{\cels\celg}\,\theta^{(N-1)}_{\celg}
\]
Hence the first equation in (\ref{equC2}) reads
\[
 \dR^\theta\pi_\gos + \frac{1}{2}\theta_{\celu_1\celu_2\gos}^{(N-3)} \wedge \Phi^{\celu_1\celu_2} + \Theta{^\celu}_{\gos\cels}\pi{_\celu}^{\cels\celg}\,\theta^{(N-1)}_{\celg} = 0
\]
We note that, since $\Theta{^\gou}_{\gog\gos} = 0$ by (\ref{equA}), the second equation in (\ref{equC2}) can also be written as
\[
 \dR^\theta\pi_\gog + \frac{1}{2}\theta_{\celu_1\celu_2\gog}^{(N-3)} \wedge \Phi^{\celu_1\celu_2} + \Theta{^\celu}_{\gog\cels}\pi{_\celu}^{\cels\celg}\,\theta^{(N-1)}_{\celg} = 0
\]
so that System (\ref{equC}) or (\ref{equC2}) can be replaced by
\begin{equation}\label{equC3}
 \dR^\theta\pi_\gou + \frac{1}{2}\theta_{\celu\celu\gou}^{(N-3)} \wedge \Phi^{\celu\celu} + \Theta{^\celu}_{\gou\cels}\pi{_\celu}^{\cels\celg}\,\theta^{(N-1)}_{\celg} = 0
\end{equation}
i.e. (\ref{ELEinsteinbrut}) with $\Lambda_0 = 0$.
We hence have recovered the constraint (\ref{constrainttoberemoved}) and Equations (\ref{ELFrobeniusdetaillee}), (\ref{ELsimple}) and (\ref{ELEinsteinbrut}), i.e. the Euler--Lagrange equations of $\mathscr{A}$.
The converse is straightforward.\hfill $\square$

%\bibliographystyle{plain}
%\bibliography{bibneo}

\end{document}